\documentclass{article}

\newcommand{\cI}{\mathcal{I}}
\newcommand{\cK}{\mathcal{K}}

\newcommand{\cP}{\mathcal{P}}

\newcommand{\cU}{\mathcal{U}}\newcommand{\cV}{\mathcal{V}}

\newcommand{\C}{\mathbb{C}}

\newcommand{\Z}{\mathbb{Z}}

\newcommand{\m}{\to}

\usepackage[vcentermath]{youngtab}
\Yboxdim{3pt} 
\newcommand{\Y}[1]{{\yng(#1)}}

\usepackage{url,graphicx,verbatim,amssymb,enumerate,stmaryrd, amsthm, mathrsfs, amsfonts, amsmath}
\usepackage{subcaption}
\usepackage[pagebackref,colorlinks,citecolor=blue,linkcolor=blue,urlcolor=blue,filecolor=blue]{hyperref}
\usepackage[dvipsnames, hyperref]{xcolor}
\usepackage{aliascnt}
\usepackage{tikz} 
\usetikzlibrary{matrix,arrows}
\usepackage{microtype}
\usepackage{multicol}
\usepackage[all]{xy}
\usepackage[hmargin=3cm,vmargin=3cm]{geometry}

\usepackage{dsfont}
\usepackage{pinlabel} 

\numberwithin{thmcounter}{section}
\newaliascnt{thmauto}{thmcounter}

\newaliascnt{Defauto}{thmcounter}

\newaliascnt{exauto}{thmcounter}

\newaliascnt{lemauto}{thmcounter}

\newaliascnt{propauto}{thmcounter}

\newaliascnt{corauto}{thmcounter}

\newaliascnt{remauto}{thmcounter}

\newaliascnt{convauto}{thmcounter}

\newtheorem{theorem}[thmauto]{Theorem}
\newtheorem{lemma}[lemauto]{Lemma}
\newtheorem{proposition}[propauto]{Proposition}
\newtheorem{corollary}[corauto]{Corollary}
\theoremstyle{definition}
\newtheorem{definition}[Defauto]{Definition}
\newtheorem{example}[exauto]{Example}
\newtheorem{remark}[remauto]{Remark}

\newtheorem*{Acknowledgments*}{Acknowledgments}

\newcommand{\mr}[1]{{\rm #1}}

\newcommand{\coker}{\mr{coker}}
\newcommand{\Hom}{\mathrm{Hom}}
\newcommand{\FI}{\mathsf{FI}}
\newcommand{\PSL}{\mathsf{PSL}}

\newcommand{\sgn}{\mathsf{sgn}}
\newcommand{\Ext}{\mathrm{Ext}}
\newcommand{\FB}{\mathsf{FB}}

\newcommand{\im}{\mathrm{im}}

\newcommand{\Tor}{\mathrm{Tor}}

\newcommand{\HH}{\mathbf H}
\newcommand{\Ind}{\mathrm{Ind}}

\newcommand{\Pun}{\mathrm{Pun}}
\newcommand{\Mod}{\mathrm{Mod}}

\DeclareMathOperator*{\colim}{\mathrm{colim}}

\title{$\FI$-hyperhomology and ordered configuration spaces}
\author{Jeremy Miller\thanks{Jeremy Miller was supported in part by NSF grants DMS-1709726 and DMS-2504473 and a Simons Foundation Travel Support for Mathematicians grant.} \, and Jennifer C. H. Wilson\thanks{Jennifer Wilson was supported in part by NSF grant DMS-1906123 and NSF CAREER grant DMS-2142709.}}

\date{\today}

\begin{document}

\maketitle

 \setcounter{tocdepth}{4}
\setcounter{secnumdepth}{4}
 
 \abstract{Using a result of Gu\`es on $\FI$-hyperhomology (a variant of a result of Gan--Li) and a semi-simplicial resolution of configuration spaces due to Kupers--Miller--Tran, we establish an explicit stable range for cohomological representation stability for configuration spaces of distinct ordered points in a manifold. Our bounds on generation degree were a factor of 5/2 better than those previously established at the time this paper was first posted.

 This updated article replaces an older version, which contained a mistake pointed out to us by Nicolas Gu\`es. We are grateful to Gu\`es for identifying the error and suggesting the resolution. 
 }
 
%\tableofcontents
 
\section{Introduction} \label{intro}

In this paper, we prove a representation stability result for ordered configuration spaces of points in a manifold. This result follows previous results of Church \cite{Ch}, Church--Ellenberg--Farb \cite{CEF}, Church--Ellenberg--Farb--Nagpal \cite{CEFN}, Miller--Wilson \cite{MillerWilson}, Church--Miller--Nagpal--Reinhold \cite{CMNR} and Bahran \cite{Cihan} to establish improved stable ranges for the cohomology of ordered configuration spaces of manifolds with integral coefficients. Bahran \cite{Cihan2024} and Gu\`es \cite{Gu} have established further improvements for manifolds satisfying additional hypotheses. 

Let $M$ be a connected manifold of dimension at least $2$. The family of $k$-point ordered configuration spaces of $M$ is not cohomologically stable in $k$; however, their cohomology groups do stabilize with respect to their natural $S_k$--actions. The standard way to define and quantify this stability is through bounds on the \emph{generation degree} and \emph{presentation degree} of the corresponding \emph{$\FI$-module} structure on the sequences of cohomology groups, which we recall below. These are equivariant analogues of the classical homological stability stable ranges in which the stabilization map is surjective, or is an isomorphism, respectively.

 Let $\FI$ denote the category of finite sets and injective maps. An \emph{$\FI$-module} is a functor from $\FI$ to the category of abelian groups, an \emph{$\FI$-chain complex} is a functor from $\FI$ to chain complexes over $\Z$, and so forth. Given an $\FI$-module $\cV$, we denote its value on a finite set $S$ by $\cV_S$. Following the notation of Church--Ellenberg--Farb \cite{CEF}, we let
$$ H_0^{\FI}: \FI\text{--}\Mod \longrightarrow \FI\text{--}\Mod $$ be the functor given by the formula $$H_0^{\FI}(V)_S :=\coker \left( \bigoplus_{T \subsetneq S} V_T \m V_S \right).$$ Here the map $V_T \m V_S$ is induced by the inclusion $T \hookrightarrow S$. Let $H_i^{\FI}$ denote the $i$th left-derived functor of $H_0^{\FI}$. These functors are known as the \emph{$\FI$-homology} functors; see Church--Ellenberg \cite{CE}. We say that  an $\FI$-module $\cV$ has \emph{generation degree $\leq g$} if $$H_0^{\FI}(\cV)_S \cong 0 \qquad \text{ for $|S|>g$.}$$ We say that $\cV$ has \emph{presentation degree $\leq r$} if $$H_0^{\FI}(\cV)_S \cong H_1^{\FI}(\cV)_S \cong 0 \qquad \text{ for $|S|>r$. }$$
Presentation degree controls the degree of generation of the terms in a 2-step resolution of $\cV$ by ``free'' $\FI$-modules; see Church--Ellenberg \cite[Proposition 4.2]{CE}.

Let $F(M)$ be the contravariant functor from $\FI$ to the category of spaces sending a finite set $S$ to the space $F_S(M)$ of embeddings of the discrete set $S$ into $M$. For the set $[k]:=\{1, 2, \ldots, k\}$, the space $F_{[k]}(M)$ is canonically identified with the \emph{ordered configuration space of $M$}, $$ F_k(M) = \{ (m_1, m_2, \ldots, m_k) \in M^k \; | \; m_i \neq m_j \text{ for all } i \neq j\}. $$ Taking integral cohomology yields $\FI$-modules $H^j(F(M))$. Our main theorem is the following.

\begin{theorem} 
\label{bounds} 
Let $M$ be a connected manifold of dimension $d = 2$. Then the $\FI$-module $H^j(F(M))$ has generation degree $\leq 4 j$ and presentation degree $\leq 4 j+1$. If $\dim M \geq 3$, then $H^j(F(M))$ has generation degree $\leq 2 j$ and presentation degree $\leq 2 j+1$.
\end{theorem}

Since the original release of this paper, but before the release of the corrected version, 
Bahran \cite[Theorem G]{Cihan2024} and Gu\`es \cite[Theorem B]{Gu} obtained refinements of \autoref{bounds}. Their results include improved ranges for more highly connected manifolds and certain surfaces. Bahran's results \cite[Example 4.12]{Cihan2024} shows the ranges in \autoref{bounds}  cannot be improved in full generality; see \autoref{CihanExample}.  

\autoref{bounds} implies that,  for $k \geq 4j$,  the group $H^j(F_k(M))$ is spanned by the $S_k$--orbit of the image of $H^j(F_{4j}(M))$ under the map induced by the fibration $F_k(M) \to F_{4j}(M)$ forgetting all but the first $4j$ points. In other words, every cohomology class in $H^j(F_k(M))$ can be represented by a sum of classes pulled back along the natural forgetful maps $F_k(M) \to F_{4j}(M)$. Moreover, for each $j$ the sequence of cohomology groups $\{H^j(F_k(M))\}_k$ is completely determined by the first $4j+1$ groups and maps between them. When $M$ has dimension at least 3, we obtain the improved result that the sequence of cohomology groups only depends on the first $2j+1$ groups and maps between them. Concretely, Church--Ellenberg \cite[Theorem C]{CE} implies the following corollary.

\begin{corollary}
Let $M$ be a connected $2$-manifold. Then for each cohomological degree $j \geq 0$, 
$$ H^j(F_k(M)) \cong \colim_{\substack{ S\subseteq \{1, 2, \ldots, k\} \\ |S| \leq 4 j+1}} H^j(F_S(M)). $$ 
If $M$ is a connected manifold of dimension $d \geq 3$, then for each cohomological degree $j \geq 0$, 
$$ H^j(F_k(M)) \cong \colim_{\substack{ S\subseteq \{1, 2, \ldots, k\} \\ |S| \leq 2 j+1}} H^j(F_S(M)). $$ 
\end{corollary} 

Another consequence is that the sequences of cohomology groups  $\{H^j(F_k(M))\}_k$ are \emph{centrally stable} in the sense of Putman \cite{Putman.Congruence}; see also Patzt \cite{Patzt.CentralStability}. This implies that, for $k \geq 4j+2$, there is an explicit presentation for the group $H^j(F_k(M))$ in terms of the groups $H^j(F_{k-1}(M))$ and $H^j(F_{k-2}(M))$, and the natural maps between them  \cite[Proposition 2.4]{CMNR}.  We therefore obtain recursive descriptions of the cohomology groups $H^j(F_k(M))$ in the stable range. 

See \cite[Theorem A]{CE} and \cite[Proposition 3.1]{CMNR} for more consequences of bounds for generation degree and presentation degree.

Another motivation for seeking improved stable ranges is the notion of \emph{secondary} stability introduced by Galatius--Kupers--Randal-Williams \cite{GKRW1, GKRW2}, which was adapted to the representation stability context in Miller--Wilson \cite{MillerWilson}. Secondary stability is a pattern in the unstable (co)homology of spaces whose (co)homology exhibits homological or representation stability. A prerequisite for establishing nontrivial secondary stability is sharp stability bounds for primary stability.

\begin{remark} This paper originated as the appendix to the preprint (arXiv Version 1)  of Miller--Wilson \cite{MillerWilson}. It was later split from that paper for reasons of length. That appendix used different methods and proved much worse bounds (exponential in cohomological degree), but provided the first explicit bounds for the generation and presentation degree of the integral cohomology of ordered configuration spaces of (possibly closed) manifolds. It was also the first paper to establish stability for configuration spaces of points in non-orientable manifolds. 
The appendix relied on algebraic results of Church--Ellenberg \cite{CE}. Since its release, Church--Miller--Nagpal--Reinhold \cite{CMNR} improved on the stable ranges, using the appendix to address non-orientable manifolds. Bahran \cite{Cihan, Cihan2024} further improved the stable ranges for high-dimensional orientable manifolds. This corrected version of the paper now relies on an algebraic result of Gu\`es \cite{Gu}, which improves the constant term in the stable range. This result of \cite{Gu} extends results of Gan--Li \cite{GanLiCongruenceSubgroups}. It applies to $\FI$-chain complexes that are not necessarily bounded from below. 
\end{remark}

\begin{remark}
Church--Miller--Nagpal--Reinhold \cite{CMNR} described two methods for proving representation stability, which they called Type A and Type B stability setups. Almost all known stability results for $\FI$-modules fall into one of those two categories broadly construed or involve directly calculating the $\FI$-modules in question. Gan--Li \cite{GanLiCongruenceSubgroups} dramatically improved the stable ranges arising from Type B setups. The methods of this paper can be used to transform most Type A setups currently in the literature into Type B setups, and should result in improved ranges in those situations. In particular, the stable ranges for pure mapping class groups (Jim\'enez Rolland \cite{RitaLinear}), homotopy groups of configuration spaces (Kupers--Miller \cite{KupersMiller.Homotopy}),  configuration spaces of non-manifolds (Tosteson \cite{Tos}), and generalized configuration spaces where collisions are allowed (Petersen \cite{Pet}) may admit improvements using the techniques of this paper. See Gu\`es \cite[Theorem A]{Gu} for subsequent work on improved stable ranges for homotopy groups of configuration spaces.
\end{remark}

\begin{remark}
 The originally published version of this paper (posted as arXiv Version 2) contains an error, which was pointed out to us by Nicolas Gu\`es. In that version, in the proof of Theorem 1.1, we incorrectly characterized the homology of the truncated complexes, and the hyperhomology spectral sequence may not collapse at the $E^2$ page as claimed.
 %Our approach to proving Theorem 1.1 did not work as written.
 Gu\`es suggested an alternative proof of the main result using ideas from his paper \cite{Gu}. 
 
% Although our approach is repairable, our main result is now superseded by more recent papers including Bahran \cite{Cihan2024} and Gu\`es \cite{Gu}.
\end{remark}

\begin{Acknowledgments*} We would like to thank Nicolas Gu\`es for pointing out the mistake in the original version of this paper and suggesting \cite[Theorem 2.2]{Gu} as a way to fix it.

We are grateful to Sarah Anderson and Rohit Nagpal for helpful conversations and to an anonymous referee for helpful feedback. 

We used ChatGPT 5.6 to proofread and give feedback on a draft of this updated article. We thank ChatGPT for drawing our attention to \cite[Example 4.12]{Cihan2024} and the related work Feichtner--Ziegler \cite{FZ}, finding  Spaltenstein \cite[Corollary 3.5]{Spaltenstein}, as well as identifying several minor errors and outdated citations. 
\end{Acknowledgments*}
 
\section{Configuration spaces of noncompact manifolds}
 
 In this section, we prove a stability theorem for the configuration spaces of a manifold $M$ in the case that $M$ is not compact. We begin by recalling the category $\FI\sharp$ which acts up to homotopy on the configuration spaces of points in a noncompact manifold. The following definition is equivalent to Church--Ellenberg--Farb \cite[Definition 4.1.1]{CEF}. 

\begin{definition} 
 Let $\FI\sharp$ denote the category whose objects are based sets $S_+ := S \sqcup \{*\}$. The morphisms are based maps with the property that the preimage of every non-basepoint element has size at most one. 
 \end{definition}
 
We will view $\FI\sharp$-modules as $\FI$-modules by restricting to the wide subcategory $\FI \subseteq \FI\sharp$ of injective morphisms. Note that the opposite category of $\FI\sharp$ is again $\FI\sharp$. 

When $M$ is a noncompact manifold, its configuration spaces admit maps $F_k(M) \to F_{k+1}(M)$ by introducing a point ``at infinity'' (see Church--Ellenberg--Farb \cite[Section 6.4]{CEF}). These maps endow the (co)homology groups with $\FI\sharp$-module structures. The following result is due to Church--Ellenberg--Farb \cite[Proposition 6.4.2]{CEF}. Miller--Wilson \cite[Section 3.1]{MillerWilson} addresses the case that the noncompact manifold $M$ is not the interior of a compact manifold with boundary. 

\begin{proposition}[{Church--Ellenberg--Farb \cite[Proposition 6.4.2]{CEF}; see also Miller--Wilson \cite[Section 3.1]{MillerWilson}}] \label{NonCompactSharp}
If $M$ is a noncompact manifold of dimension $d \geq 2$, then the assignment $S_{+} \mapsto F_S(M)$ defines a functor from $\FI\sharp$ to the homotopy category of spaces. In particular, post-composing with the $j$th homology or cohomology group functor gives $\FI\sharp$-modules which we denote by $H_j(F(M))$ and $H^j(F(M))$, respectively.
\end{proposition}

The following result appears in Miller--Wilson \cite[Theorem 3.12 and Theorem 3.27]{MillerWilson}. It follows easily from Church--Ellenberg--Farb \cite[Theorem 6.4.3]{CEF} if $M$ is orientable and of finite type. 

\begin{theorem}[{\cite[Theorems 3.12 and 3.27]{MillerWilson}, \cite[Theorem 6.4.3]{CEF}}] \label{ConfigSpaceRepStable} Let $M$ be a noncompact connected manifold of dimension $d \geq 2$. If $d=2$, then $H_j(F(M))$ has generation degree $\leq 2j$. If $d \geq 3$, then $H_j(F(M))$ has generation degree $\leq j$. 
\end{theorem}

\begin{definition}\label{DefnInduced} 
Let $\FB$ denote the category of finite sets and bijections. Let $\Mod_{\FB}$ denote the category of functors from $\FB$ to abelian groups, so $\Mod_{\FB}$ is equivalent to the category of sequences of integral $S_k$--representations. Let $\cI: \Mod_{\FB} \m \Mod_{\FI}$ be the left adjoint of the forgetful functor. Concretely, for an $\FB$-module $\cU$, the $\FI$-module $\cI(\cU)$ is given by the following formula.  See Church--Ellenberg--Farb \cite[Definition 2.2.2]{CEF}, where $\cI(\cU)$ is denoted $M(\cU)$.   
\begin{equation} \label{InducedModuleFormula}
\cI(\cU)_A := \bigoplus_{B\subseteq A} \cU_{B}. 
\end{equation} 
Observe 
$$ \cI(\cU)_A  \quad \cong \quad \bigoplus_{ 0 \leq b \leq |A|} \Ind_{S_b \times S_{|A|-b}}^{S_{|A|}} \cU_{[b]} \boxtimes \mathrm{Trv}_{|A|-b} \qquad (\mathrm{Trv}_{|A|-b} \text{ the trivial representation).} 
$$
$\FI$-modules of the form $\cI(\cU)$ are called \emph{induced} $\FI$-modules. 
\end{definition} 

Using the functor $\cI$, Church--Ellenberg--Farb \cite[Theorem 4.1.5]{CEF} gave another characterization of $\FI\sharp$-modules.

\begin{theorem}[{\cite[Theorem 4.1.5]{CEF}}] \label{FIsharp} The functor $\cI$ factors through the inclusion $\Mod_{\FI\sharp} \m \Mod_{\FI}$ and induces an equivalence of categories $\Mod_{\FB} \m \Mod_{\FI\sharp}$. 

\end{theorem}

Church--Ellenberg--Farb \cite[Remark 2.3.8]{CEF} showed that, if we view $\Mod_{\FB}$ as the codomain of $H_0^{\FI}$, then $H_0^{\FI}$ and $\cI$ are quasi-inverse on $\FI\sharp$-modules. 

For an $\FI$-module or  $\FB$-module $\cU$, we say $\deg \cU \leq d$ if $\cU_k \cong 0$ for $k>d$. In particular, the generation degree of $\cI(\cU)$ is equal to $\deg \cU$. 
 
 Given an $\FI\sharp$-module $\cV$, since $\FI\sharp$ is self-dual, we can form $\FI\sharp$-modules $\Ext^i_\Z(\cV,\Z)$ by composing with the $\Ext$ functors pointwise. Since the opposite category of $\FB$ is $\FB$, we can also perform this construction with $\FB$-modules to obtain $\FB$-modules. 
 
 \begin{lemma} \label{ExtLemma}
 Let $\cV$ be an $\FI\sharp$-module with generation degree $\leq g$. Then the $\FI\sharp$-modules $\Ext^i_\Z(\cV,\Z)$ have generation degree $\leq g$. 
  \end{lemma}
 
 \begin{proof}
Since the direct sum of exact chain complexes of free abelian groups is an exact chain complex of free abelian groups, it follows from the description of $\cI$ in Equation (\ref{InducedModuleFormula}) (quoting \cite[Definition 2.2.2]{CEF}) that $\cI$ takes pointwise free resolutions to pointwise free resolutions. Similarly, since $\Hom$ commutes with finite direct sums, it follows that
$\cI$ commutes with pointwise $\Hom$ functors. Thus $\cI$ commutes with pointwise $\Ext$ functors and we have isomorphisms 
 
$$ \Ext^i_\Z(\cV,\Z) \cong \Ext^i_\Z(\cI(H_0^{\FI}(\cV)  ),\Z) \cong \cI \big( \Ext^i_\Z(H_0^{\FI}(\cV)  ,\Z) \big).$$ To check that these isomorphisms of abelian groups assemble to form an isomorphism of $\FI\sharp$-modules, it is helpful to use functorial free resolutions like the bar resolution to compute $\Ext$. This isomorphism shows that the generation degree of $\Ext^i_\Z(\cV,\Z)$ is bounded by $\deg H_0^{\FI}(\cV)$ which is the generation degree of $\cV$.
 \end{proof}

 This gives the following corollary.

\begin{corollary} \label{RepStabilityCohomology} Let $M$ be a noncompact connected manifold of dimension $d$. If $d=2$, then $H^j(F(M))$ has generation degree $\leq 2j$. If $d \geq 3$, then $H^j(F(M))$ has generation degree $\leq j$.
\end{corollary}

\begin{proof}
By the universal coefficient theorem, there is a short exact sequence of $\FI\sharp$-modules $$0 \m \Ext^1_\Z(H_{j-1}(F(M),\Z)) \m H^j(F(M)) \m \Ext^0_\Z(H_j(F(M)),\Z) \m 0.$$ This gives an exact sequence $$H_0^{\FI}\left( \Ext^1_\Z(H_{j-1}(F(M),\Z) \right) \m H_0^{\FI}\left(H^j(F(M)) \right) \m H_0^{\FI}\left( \Ext^0_\Z(H_j(F(M)),\Z) \right).$$ The claim now follows by \autoref{ConfigSpaceRepStable} and \autoref{ExtLemma}. 
\end{proof}
\begin{remark} We remark that, when $M$ has finite type, \autoref{RepStabilityCohomology}  is a consequence of the result of Church--Ellenberg--Farb \cite[proof of Theorem 4.1.7]{CEF} that an $\FI\sharp$-module $\cV$ is finitely generated in degree $\leq g$ if and only if $\cV_k$ is generated by $O(k^g)$ elements. We can then use \autoref{ConfigSpaceRepStable} and the universal coefficient theorem to bound the growth in the number of generators of the abelian group $H^j(F_k(M))$, and \autoref{RepStabilityCohomology} follows. 
 \end{remark} 

Higher $\FI$-homology groups of the cohomology groups of configuration spaces of noncompact manifolds vanish uniformly by the following theorem.

\begin{theorem}[{\cite[Lemma 2.3]{CE}}] \label{FIsharpAcyclic} Let $\cV$ be an $\FI\sharp$-module. Then $H^{\FI}_i(\cV)=0$ for all $i > 0$.
\end{theorem}

\section{The puncture resolution}

We now recall the \emph{puncture resolution},   which will allow us to leverage results for noncompact manifolds to prove results for compact manifolds. We use a discrete version of the semi-simplicial space from Randal--Williams \cite{RW}. Our construction is a special case of the augmented semi-simplicial space $\tilde{\mathcal{W}}_{\bullet}(\lambda)$ of Kupers--Miller--Tran \cite[Definition 5.3]{KMT}; the following definition is the case that $\lambda$ is the partition $1^{|S|-2}2$. 

\begin{figure}[h!]
\begin{center} 
\includegraphics[scale=.3]{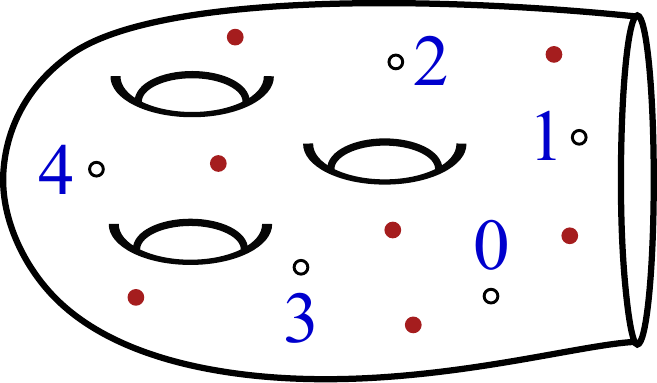}
\end{center}
\caption{An element of $\Pun_4(F_7(M))/S_7$. The $(4+1)$ ordered punctures are labelled in blue, and the 7 particles in the unordered configuration of the punctured surface are marked in red.}
\label{PunResolution.eps}
\end{figure}

\begin{definition} Fix a finite set $S$. We define the \emph{puncture resolution} of the configuration space $F_S(M)$ as follows. For $p \geq -1$, denote the topological spaces
\begin{align*} 
\Pun_p(F_S(M)) &= \bigsqcup_ {(x_0,\ldots, x_p) \in F_{\{0,\ldots,p \} }(M)} F_S(M-\{x_0,\ldots, x_p\}). 
\end{align*} 
\end{definition}
In particular, $\Pun_{-1}(F_S(M)) = F_S(M)$.  See Figure \ref{PunResolution.eps} for an illustration.
For fixed $S$, the spaces $\Pun_p(F_S(M))$ assemble to form an augmented semi-simplicial space. The $i$th face map is induced by the inclusion $$M-\{x_0,\ldots, x_p\} \hookrightarrow M-\{x_0,\ldots, \hat x_i, \ldots, x_p\}.$$ 
Given an injection $S \hookrightarrow T$, composition of embeddings gives a map of augmented semi-simplicial spaces $\Pun_\bullet(F_T(M)) \m \Pun_\bullet(F_S(M))$. In this way, the puncture resolution assembles to form an augmented semi-simplicial co-$\FI$-space $\Pun_\bullet(F(M)) $. Thus the cohomology groups of the space of $p$-simplices of the puncture resolution assemble to form an $\FI$-module and the face maps are maps of $\FI$-modules. In particular, the action of the symmetric group $S_k$ on $F_k(M)$ extends to an action on $\Pun_\bullet(F_k(M))$. 

Kupers--Miller--Tran prove the following as the special case of \cite[Proposition 5.8]{KMT} where (in their notation) $\lambda$ is the partition $1^{k-2}2$. 
 
\begin{theorem}[{\cite[Proposition 5.8]{KMT}}] \label{OscarPuncture} \label{punctureResolution} Let $M$ be a connected manifold of dimension at least $2$. The augmentation map $$||\Pun_\bullet(F_k (M))|| \m F_k(M)$$ is a weak equivalence.
\end{theorem}
 
% We now use \autoref{OscarPuncture} to prove the following.

%\begin{proposition} \label{punctureResolution} Let $M$ be a manifold. The augmentation map $||\Pun_\bullet(F_k (M))|| \m F_k(M)$ is a weak equivalence. %\end{proposition}

%\begin{proof} We proceed as in the proof of Miller--Wilson \cite[Proposition 3.8]{MillerWilson}. It suffices to show that  the homotopy fiber of the map $||\Pun_\bullet(F_k (M))|| \m F_k(M)$ is contractible. If we use the standard path space construction of homotopy fibers, then this homotopy fiber is homeomorphic to the homotopy fiber of $||\Pun_\bullet(F_k (M))||/S_k \m F_k(M)/S_k$, which is contractible by \autoref{OscarPuncture}.
%
%\end{proof}

\section{Proof of \autoref{bounds}}
 
 In this section, we prove the main theorem. We begin by recalling the definition of $\FI$-hyperhomology.
 
 \begin{definition}
 Given an $\FI$-chain complex $\cV_*$, we denote the $i$th hyperhomology group associated to the chain complex $\cV_*$ and the functor $H_0^{\FI}$ by $\HH_i^{\FI}(\cV_*)$ and refer to it as \emph{$\FI$-hyperhomology}.
 \end{definition}
 
 A key component of our proof is the following result of Gu\`es \cite[Theorem 2.2 (2)]{Gu}; also see {Gan--Li \cite[Theorem 5]{GanLiCongruenceSubgroups}}. Note that Gu\`es' result is stated in terms of homotopy groups of spectra but in the case of Eilenberg--MacLane spectra, this is equivalent to homology groups of chain complexes. We use the term \emph{ homologically graded chain complex} to mean one whose differentials decrease homological degree by $1$.

 \begin{theorem}[Gu\`es {\cite[Theorem 2.2 (2)]{Gu}}] \label{GanLi5}
 Let $\mathcal{V_*}$ be a homologically graded $\FI$--chain complex. Suppose for all $i \in \Z$ that $$\deg {\bf H}^{\FI}_i(\cV_*) \leq -a i + b $$ with $a>0$. Then, for all  $i \in \Z$, we have that 
 \begin{align*}
 \deg H_0^{\FI}(H_{-i}(\cV_*)) &\leq 2ai+2b, \\ 
  \deg H_1^{\FI}(H_{-i}(\cV_*)) &\leq 2ai+2b+1.
 \end{align*}
 \end{theorem} 

Church--Ellenberg \cite[Definition 5.9 and Proposition 5.10]{CE} prove the following result. See also Church--Ellenberg--Farb--Nagpal \cite[Definition 2.19 and Remark 2.21]{CEFN}.

 \begin{lemma}[Church--Ellenberg {\cite{CE}}] \label{FIHomologyComplex} Let $\cU$ be an $\FI$-module. Let $\sgn_p$ denote the rank-1 sign representation of the symmetric group $S_p$.  Then the $\FI$-homology groups $H_p^{\FI}(\cU)$ of $\cU$ are computed by the following $\FI$-chain complex $C_*^{\FI}(\cU)$.
 $$ C_p^{\FI}(\cU)_S := \bigoplus_{\substack{\FI \;\mathrm{ morphisms} \\ f: [p] \to S}} \cU_{S \setminus \im(f)} \otimes_{\Z[S_p]} \sgn_p .$$ The differentials are induced by the alternating sums $\sum_{i=1}^p (-1)^{i-1} d_i$ of the face maps 
 $$ d_i:  \bigoplus_{\substack{\FI \;\mathrm{ morphisms} \\ f: [p] \to S}} \cU_{S \setminus \im(f) }  \; \longrightarrow \bigoplus_{\substack{\FI \;\mathrm{ morphisms} \\ g: [p-1] \to S}} \cU_{S \setminus \im(g)} $$ 
arising from the natural semi-simplicial group structure   $\bigoplus_{ f: \bullet \to S} \cU_{S \setminus \im(f) }.$
 \end{lemma}

 Observe 
 \begin{align*} C_p^{\FI}(\cU)_S & \quad \cong \quad \Ind_{S_{k-p} \times S_{p}}^{S_k} \mathcal{U}_{k-p} \boxtimes \sgn_p 
 \end{align*}

We note that, from our description of the functor $C_p^{\FI}(-)_S$, that it (and hence $H_p^{\FI}(-)_S$) commute with arbitrary products.

 \begin{lemma} 
Let $\mathcal U_*$ be an $\FI$-chain complex, which need not be bounded above or below.  The total complex associated to the double complex $C^{\FI}_*(\mathcal U_*)$ computes the $\FI$-hyperhomology ${\bf H}^{\FI}_i(\cU_*)$.
 \end{lemma} 

 \begin{proof}  Let $\cK$ denote the $(\FI^{op} \times \FI)$-module 
 $$\cK(S,T) = \left\{ \begin{array}{ll} \Z\left[\text{Bijections}(S,T)\right], & |S|=|T|, \\  0, & |S| \neq |T| \end{array}\right.$$ 

 Church--Ellenberg \cite[Proposition 5.3]{CE} show that there are natural isomorphisms of $\FI$-modules
 $$ H_0^{\FI}(\cU) \cong \cU \otimes_{\FI} \cK \qquad \qquad H_i^{\FI}(\cU) \cong \Tor_{i}^{\FI} (\cU, \cK).$$ 
 where $\cU\otimes_{\FI} \cK$ denotes the tensor product of functors in the sense of (for example) Mac Lane \cite[Section IX.6]{MacLane}.  Church--Ellenberg \cite[Definition 5.9 \& Proposition 5.10]{CE} construct a resolution of $(\FI^{op} \times \FI)$-bimodules $C_{*} \to \cK$ satisfying the following.
 \begin{enumerate}
     \item[(i)] By its definition, for fixed $T$, the $\FI^{op}$-module $C_p(-,T)$ vanishes unless $0 \leq p \leq |T|$. 
     \item[(ii)] For each fixed $T$ and $p$, the $\FI^{op}$-module $C_p(-, T)$ is a projective $\FI^{op}$-module. %In particular, for fixed $T$, the resolution of $\FI^{op}$-modules $C_*(-, T) \to \cK(-, T)$ is a bounded and projective. 
     \item[(iii)] Given an $\FI$-module $\cU$, there is an isomorphism of $\FI$-chain complexes
     $$C_*^{\FI}(\cU) \quad \cong \quad \cU \otimes_{\FI} C_*.$$
 \end{enumerate}
 They establish properties (ii) and (iii) to prove the complex $C_*^{\FI}(\cU)$ computes $H_*^{\FI}(\cU)$, which we quote in \autoref{FIHomologyComplex} above. 

 Let $\cU_*$ be an $\FI$-chain complex. If $\cU_*$ is unbounded, it is not clear that the total complex of a naively constructed Cartan--Eilenberg resolution of $\cU_*$ will be quasi-isomorphic to $\cU_*$, nor that the associated hyperhomology spectral sequences will converge. For the present proof we require this quasi-isometry property, so we seek instead a resolution $\cP_{*} \to \cU_*$ that is \emph{$K$-projective} in the sense of Spaltenstein \cite[Definition p. 124]{Spaltenstein}. 
 
 Spaltenstein \cite[Corollary 3.5]{Spaltenstein} ensures that $K$-projective resolutions exist for chain complexes in the category of $\FI$-modules over $\Z$.  Spaltenstein's result states existence holds for chain complexes in a category as long as this category (i) has enough projectives; (ii) has filtered colimits; and (iii) $\varinjlim$ is exact. Their construction of the $K$-projective resolution is term-wise projective. We can verify that these criteria hold for $\FI$-modules. The category has enough projectives; see Church--Ellenberg--Farb \cite[Definition 2.2.3 and Remark 2.2.A]{CEF}. For the last two properties, observe that filtered colimits in $\FI$-modules correspond object-wise to filtered colimits in abelian groups, a category where both properties hold.

So, for our fixed $\FI$-chain complex $\cU_*$, let $\cP_{*}$ be $K$-projective resolution whose terms are projective $\FI$-modules. We will study the double complex $\cP_{*} \otimes_{\FI} C_*$. 

We will first compare the total complex of $\cP_{*} \otimes_{\FI} C_*$ to the complex $\cP_{*} \otimes_{\FI} \cK$ that computes the $\FI$-hyperhomology ${\bf H}^{\FI}_i(\cU_*)$. We note that, although the complex $\cP_{*}$ is term-wise projective, we do not automatically obtain a  quasi-isomorphism between $\cP_{*} \otimes_{\FI} C_*$ and $\cP_{*} \otimes_{\FI} \cK$ because
the complex $\cP_{*}$ may be unbounded. If, however, we fix an object $T$ of $\FI$, we can use the boundedness property (Property (i)) of the complex  $C_*(-,T)$ to avoid this issue. 

Observe that, for fixed $p$, there are only finitely many $i,j$ such that $i+j=p$ and  $\cP_{i} \otimes_{\FI} C_j(-,T)$ is nonzero, without any boundedness assumptions on the complexes $\cU_*$ or $\cP_*$. Thus the spectral sequence of the double complex  $\cP_{*} \otimes_{\FI} C_*(-, T)$ converges, and we deduce that its total complex is quasi-isomorphic to $\cP_{*} \otimes_{\FI} \cK(-, T)$.

On the other hand, we will compare the total complexes of $\cP_{*} \otimes_{\FI} C_*$ and $\cU_* \otimes_{\FI} C_* \cong C_*^{\FI}(\cU_*)$. 

Fix an object $T$ of $\FI$. We will establish a quasi-isomorphism between the total complexes of the double complexes $\cP_{*} \otimes_{\FI} C_*(-,T)$ and $C_*^{\FI}(\cU_*)_T$.  Again the boundedness property (Property (i)) of the complex  $C_*(-,T)$ ensures that the spectral sequence associated to the double complex $\cP_{*} \otimes_{\FI} C_*(-,T)$ converges. Hence, using the projectivity of $C_*(-,T)$ (Property (ii)), we obtain the desired quasi-isomorphism, and complete the proof. 
 \end{proof}

\begin{lemma} \label{greaterthanlemma}
Let $\mathcal V_* \to \mathcal W_*$ be a map of $\FI$-chain complexes such that $H_j(\mathcal V_*) \to H_j(\mathcal W_*)$ is an isomorphism for $j \geq d$. Then ${\bf H}^{\FI}_i(\mathcal V_*)_k \to {\bf H}^{\FI}_i(\mathcal W_*)_k$ is an isomorphism for $i \geq d+k+1$.
\end{lemma}

\begin{proof}
%Given an $\FI$-module $\mathcal U$, let $$C^{\FI}_i(\mathcal U)_k = \Ind_{S_{k-i} \times S_{i}}^{S_k} \mathcal{U}_{k-i} \boxtimes \sgn_i$$ with $\sgn_i$ the sign representation. See Church--Ellenberg \cite[Definition 5.9 and Proposition 5.10]{CE} for a description of a differential $ C^{\FI}_i(\mathcal U)_k  \to C^{\FI}_{i-1}(\mathcal U)_k  $ with the property that $$H_i(C^{\FI}_*(\mathcal U) )={\bf H}_i^{\FI}(\mathcal U).$$ For formal reasons, if $\mathcal U_*$ is an $\FI$-chain complex, the homology of the total complex associated to the double complex $C^{\FI}_*(\mathcal U_*)$ computes the $\FI$-hyperhomology.

Consider a map of $\FI$-chain complexes $\mathcal V_* \to \mathcal W_*$ such that $H_j(\mathcal V_*) \to H_j(\mathcal W_*)$ is an isomorphism for $j \geq d$. Let $\mathcal Z_*$ denote the mapping cone. Note that $H_j(\mathcal Z_*) \cong 0$ for $j \geq d+1$. Thus, $H_j(C_i^{\FI}(\mathcal Z_*)) \cong 0$ for $j \geq d+1$. Since $C_i^{\FI}(\mathcal Z_*)_k =0$ for $i>k$, $H_j(C_i^{\FI}(\mathcal Z_*))_k \cong 0$ for $i>k$. Consider the spectral sequence: $$E^0_{p,q}=C^{\FI}_{p}( \mathcal Z_q)_k \implies {\bf H}^{\FI}_{p+q}(\mathcal Z_*)_k.$$ 
Note that $E^1_{p,q} \cong 0$ for $q \geq d+1$ or $p \geq k+1$. Thus, $E^1_{p,q} \cong 0$ for $p+q \geq d+k+1$ and hence $${\bf H}^{\FI}_{i}(\mathcal Z_*)_k \cong 0 \text{ for } i \geq d+k+1.$$ The claim now follows from the long exact sequence in $\FI$-hyperhomology.
\end{proof}

 The following result is well-known; see for example Bendersky--Gitler \cite[Proposition 1.2]{BenderskyGitler}. 
 
 \begin{proposition}[{E.g.,  Bendersky--Gitler \cite[Proposition 1.2]{BenderskyGitler}}] \label{BG}
 Let $X_\bullet$ be a semi-simplicial space. Let $X^{*,*}$ be the cohomologically graded double complex with $X^{p,q}=C^q(X_p)$ with one differential given by the differential on singular cochains and the other given by the alternating sum of the face maps. Let $X^*$ be the total complex of $X^{*,*}$. Then $H^i(X^*) \cong H^i(||X_\bullet||)$.

 \end{proposition}

We now prove a vanishing line for the $\FI$-hyperhomology of configuration spaces. Since \autoref{GanLi5} requires homologically graded chain complexes, we cannot apply it to $C^*(F(M))$. Instead we consider the homologically graded $\FI$-chain complex $C^{-*}(F(M))$ which vanishes in positive homological degrees. 

\begin{proposition} \label{VanishFI}
Let $M$ be a connected manifold with $\dim M \geq 2$. Then $\deg {\bf H}^{\FI}_i(C^{-*}(F(M))) \leq -2 i$. If $\dim M \geq 3$, then $\deg {\bf H}^{\FI}_i(C^{-*}(F(M))) \ \leq -i$.
\end{proposition}

\begin{proof}
 Let $r \geq 0$. Observe that the spaces $\Pun_r(F_k(M)) $ assemble to form a co-$\FI$-space which we denote by $\Pun_r(F(M)) $. Thus, $C^{-*}(\Pun_r(F_k(M))) $ is an $\FI$-chain complex. It is isomorphic to a direct product of $\FI$-chain complexes of the form   $C^{-*}(F(N)) $ with $N$ equal to $M$ with a nonzero number of points removed. In particular, $N$ is noncompact. \autoref{RepStabilityCohomology} now implies $$\deg H_0^{\FI} \left( H_j( C^{-*}(F(N))) \right) = \deg H_0^{\FI} \left( H^{-j}( F(N)) \right) \leq - \alpha j  \quad \text{ for } \quad  \alpha = \left\{ \begin{array}{ll} 2, & \dim(M)=2 \\ 1, & \dim(M) \geq 3. \end{array} \right.$$ 

By \autoref{FIsharpAcyclic}, $$H_i^{\FI} \left( H_j( C^{-*}(F(N))) \right) = H_i^{\FI} \left( H^{-j}( F(N)) \right)=0 \qquad \text{ for $i>0$.}$$
Since homology commutes with direct products (Weibel \cite[Exercise 2.1.2]{Weibel}), it follows that  $$ { H}_i^{\FI} \left( H^{-j}(\Pun_r(F(M))) \right)_k  \cong 0 \qquad \text{ if $i>0$ or $k > -\alpha j$.}$$
For each fixed $r \geq 0$, one of the two hyperhomology spectral sequences associated with the chain complex $C^{-*}(\Pun_{r}(F(M)))$  satisfies $$E^2_{p,q}(k) =  { H}^{\FI}_p \left( H^{-q}(\Pun_r(F(M))) \right)_k  \implies {\bf H}^{\FI}_{p+q}\left( C^{-*}(\Pun_r(F(M)) \right)_k.$$ 
It collapses at $E^2$ and we see that 
\begin{equation} \label{HyperhomologyPun} {\bf H}_{i}^{\FI} \left( C^{-*}(\Pun_r(F(M))) \right)_k \cong { H}_0^{\FI} \left( H^{-i}(\Pun_r(F(M))) \right)_k
\end{equation} which vanishes if $k > -\alpha i$.

By \autoref{punctureResolution} and \autoref{BG}, $C^{*}(F(M))$ is quasi-isomorphic to the total complex of the cohomologically graded $\FI$-double complex with $(p,q)$-entry $ C^q \left( \Pun_p(F(M)) \right)$. Let $$\mathcal{C}_{p,q}=C^{-q} \left( \Pun_{-p}(F(M)) \right)$$ be the $\FI$-double complex obtained by negating the indices and let $\mathcal{C}_*$ denote the total complex. Since $C^{-*}(F(M))$ is equivalent to $\mathcal{C}_*$, our goal is to show ${\bf H}^{\FI}_i(\mathcal{C}_*)_k$ vanishes for $k > -i \alpha$. Let $\mathcal F^r_*$ be the sub $\FI$-chain complex of $\mathcal C_*$ associated to the portion of the double complex $\mathcal{C}_{p,q}$ with $p \leq r$. Note that  $\mathcal F^0_* = \mathcal C_*$ and $\mathcal F^{r+1}_* \supseteq \mathcal F^{r}_*$. 

The map $H_j(\mathcal C_*) \to H_j(\mathcal C_*/\mathcal F^r_*)$ is an isomorphism for $j \geq r+2$ since $\mathcal C_j \to \mathcal C_j/\mathcal F^r_j$ is an isomorphism for $j \geq r+1$. Thus, for $r$ sufficiently negative (specifically $r \leq i-k-3$), \autoref{greaterthanlemma} implies ${\bf H}^{\FI}_i(\mathcal C)_k \to {\bf H}^{\FI}_i(\mathcal C_*/\mathcal F^r_*)_k$ is an isomorphism. In particular, it suffices to show ${\bf H}^{\FI}_i(\mathcal C_*/\mathcal F^r_*)_k \cong 0$ for $k > -\alpha i$ and all $r$.

We will prove this by induction on $r$. For $r \geq 0$, the complex $\mathcal C_*/\mathcal F^r_*$ vanishes and the result is vacuous. Fix $r<0$ and assume we have proven the claim for $r+1$. We have a short exact sequence of $\FI$-chain complexes
$$ 0 \; \longrightarrow \; \mathcal F^{r+1}_* /\mathcal F^{r}_*   \; \longrightarrow \; \mathcal C_*/\mathcal F^{r}_* \; \longrightarrow \; \mathcal C_*/\mathcal F^{r+1}_* \; \longrightarrow \;0.$$ Observe that  $$\mathcal F^{r+1}_* /\mathcal F^{r}_* \cong C^{-*+r+1}(\Pun_{-r-1}(F(M))).$$ Thus  $${\bf H}^{\FI}_i(\mathcal F^{r+1}_* /\mathcal F^{r}_* )_k \cong {\bf H}_{i-r-1}^{\FI}\left(C^{-*}(\Pun_{-r-1}(F(M))) \right)_k$$ which vanishes for $k> -\alpha (i-r-1)=-\alpha i +\alpha (r+1)$ by our analysis of Equation (\ref{HyperhomologyPun}). Since $r+1 \leq 0$ and $\alpha > 0$, ${\bf H}^{\FI}_i(\mathcal F^{r+1}_* /\mathcal F^{r}_* ) \cong 0$ for $k> -\alpha i$. It follows from the long exact sequence in $\FI$-hyperhomology that ${\bf H}^{\FI}_i(\mathcal C_*/\mathcal F^r_*)_k \cong 0$ for $k > -\alpha i$.
\end{proof}

\autoref{bounds} now follows from \autoref{GanLi5} and \autoref{VanishFI}.
\begin{example} \label{CihanExample} We remark that the bounds on generation and presentation degree in \autoref{bounds} are sharp in the case of $M=S^2$ and cohomological degree $1$, as observed by Bahran \cite[Examples 1.13 and 4.12]{Cihan2024}. Bahran's example,  and the integral cohomology computation of $F_n(S^2)$ due to Feichtner--Ziegler \cite[Theorem 2.4]{FZ}, imply that the integral $\FI$-module $H^1(F(S^2))$ has a partial resolution by induced $\FI$-modules 
$$\cI\left(V_{\Y{2,2,1}}\right) \longrightarrow \cI\left(V_{\Y{2,2}}\right) \longrightarrow H^1(F(S^2)) \longrightarrow 0 $$
where $V_{\Y{2,2}}$ denotes the $\FB$-module that has a copy of the $S_4$-representation $V_{\Y{2,2}} \cong \Z^2$ in degree $4$ and zero elsewhere, and  $V_{\Y{2,2,1}}$ denotes the $\FB$-module with a copy of the $S_5$-representation $V_{\Y{2,2,1}} \cong \Z^5$ in degree $5$ and zero elsewhere. 

\begin{comment}
In fact, the calculation of the integral cohomology algebra of $F_n(S^2)$ due to Feichtner--Ziegler \cite[Theorem 2.4]{FZ} shows that $H^1(F_n(S^2))$ is free abelian for all $n$ with
$$H_0^{\FI}\left(H^1\left(F\left(S^2\right)\right)\right)_n \cong \left\{ \begin{array}{ll} \Z^2, & n=4 \\ 0, & \text{otherwise}.\end{array} \right.  $$ 
so the $\FI$-module $H^1(F(S^2))$ has generation degree exactly 4. Moreover, $H^1(F_n(S^2))$ is free abelian with
$$ \mathrm{rank}\left(H_1\left(F_5\left(S^2\right)\right)\right) = {5-1 \choose 2}-1 = 5 \quad < \quad 10 = \mathrm{rank}\left( \Ind_{S_4 \times S_1}^{S_5} \Z^2 \right)$$ and we see that $H_1(F(S^2))$ is not a induced $\FI$-module (\autoref{DefnInduced}) and so is not projective; see \cite[Definition 2.2.2 and Remark 2.2.A]{CEF}. \autoref{bounds} then implies that $H_1(F(S^2))$ has presentation degree exactly 5. It has a partial resolution of the form
$$\cI(\Z^5[5]) \longrightarrow \cI(\Z^2[4]) \longrightarrow H^1(F(S^2)) \longrightarrow 0 $$
where $\Z^2[4]$ denotes a $\FB$-module with underlying abelian group $\Z^2$ in degree 4 and zero elsewhere; $\Z^5[5]$ is analogous. 
%where the $\FB$-modules $\cU^0$ and $\cU^1$ have underlying abelian groups $$\cU^0 \cong \left\{ \begin{array}{ll} \Z^2, & n=4 \\ 0, & \text{otherwise}.\end{array} \right. \qquad \text{and} \qquad \cU^1 \cong \left\{ \begin{array}{ll} \Z^5, & n=5 \\ 0, & \text{otherwise}.\end{array} \right.$$
\end{comment} 
\end{example}

{\footnotesize
\bibliographystyle{amsalpha}
\bibliography{hyper}
}

\end{document}